\documentclass[12pt]{article}
\usepackage{cite}
\usepackage{mathrsfs}
\usepackage{amssymb,amsfonts,amsmath}
\usepackage{enumerate}
\usepackage{indentfirst}
\usepackage{latexsym,bm}
\usepackage[noend]{algpseudocode}
\usepackage{algorithmicx,algorithm}
\usepackage{algorithm}
\usepackage{makeidx}
\usepackage{fancybox}
\usepackage{color}
\usepackage{multicol}
\usepackage[latin1]{inputenc}
\usepackage{graphicx}
\usepackage{epstopdf}
\usepackage{pstricks,pst-node,pst-text,pst-3d}
\usepackage{tikz}
\usepackage{amsfonts,amsmath,amssymb}
\usepackage{tikz}
\usetikzlibrary{shapes,chains,scopes,positioning,backgrounds,fit,shadows,calc,arrows.meta}
\usepackage{caption}
\newtheorem{thm}{Theorem}[section]

\newtheorem{lem}[thm]{Lemma}

\newtheorem{pro}[thm]{Proposition}

\newenvironment{pf}{{\noindent \it \bf Proof:}}{{\hfill$\Box$}\\}

\def\qed{\hfill \nopagebreak\rule{5pt}{8pt}}

\begin{document}

\title{\bf Directed Steiner path packing and directed path
connectivity}
\author{
Yuefang Sun\\
School of Mathematics and Statistics, Ningbo University,\\
Ningbo 315211, P. R. China\\ 
Email address: sunyuefang@nbu.edu.cn}
\maketitle

\begin{abstract}
For a digraph $D=(V(D), A(D))$, and a set $S\subseteq V(D)$ with $r\in S$ and $|S|\geq 2$, a directed $(S, r)$-Steiner path or, simply, an $(S, r)$-path is a directed path $P$ started at $r$ with $S\subseteq V(P)$. Two $(S, r)$-paths are said to be arc-disjoint if they have no common arc. Two arc-disjoint $(S, r)$-paths are said to be internally disjoint if the set of common vertices of them is exactly $S$. Let $\kappa^p_{S,r}(D)$ (resp. $\lambda^p_{S,r}(D)$) be the maximum number of internally disjoint (resp. arc-disjoint) $(S, r)$-paths in $D$. The directed path $k$-connectivity of $D$ is defined as
$$\kappa^p_k(D)= \min \{\kappa^p_{S,r}(D)\mid S\subseteq V(D), |S|=k, r\in
S\}.$$ Similarly, the directed path $k$-arc-connectivity of $D$ is defined as
$$\lambda^p_k(D)= \min \{\lambda^p_{S,r}(D)\mid S\subseteq V(D), |S|=k, r\in
S\}.$$ The directed path $k$-connectivity  and directed path $k$-arc-connectivity are also called directed path connectivity which
extends the path connectivity on undirected graphs to
directed graphs and could be seen as a generalization of classical connectivity of digraphs.

In this paper, we study the complexity for $\kappa^p_{S,r}(D)$ and $\lambda^p_{S,r}(D)$.
When both $k\geq 3, \ell\geq 2$ are fixed integers, we show that the problem of deciding whether $\kappa^p_{S,r}(D) \geq \ell$ for an Eulerian digraph $D$ is NP-complete, where $r\in S\subseteq V(D)$ and $|S|=k$. However, when we consider the class of symmetric digraphs, the problem becomes polynomial-time solvable. We also show that the problem of deciding whether $\lambda^p_{S,r}(D) \geq \ell$ for a given digraph $D$ is NP-complete, where $r\in S\subseteq V(D)$ and $|S|=k$.

We also give two upper bounds for $\kappa^p_k(D)$ in terms of classical connectivity $\kappa(D)$ and $\lambda(D)$  of digraphs, and give sharp lower and upper bounds for the Nordhaus-Gaddum type relations of the parameter $\lambda^p_k(D)$.

\vspace{0.3cm}

\vspace{0.3cm}{\bf AMS subject
classification (2020)}: 05C20, 05C38, 05C45, 05C70, 05C85, 68Q25, 68R10.

\end{abstract}

\section{Introduction}

For a graph $G=(V,E)$ and a (terminal) set $S\subseteq V$ of at
least two vertices, an {\em $S$-Steiner tree} or, simply, an {\em
$S$-tree} is a tree $T$ of $G$ with $S\subseteq V(T)$. Two $S$-trees $T_1$ and $T_2$ are said to be {\em
edge-disjoint} if $E(T_1)\cap E(T_2)=\emptyset$. Two edge-disjoint
$S$-trees $T_1$ and $T_2$ are said to be {\em internally disjoint}
if $V(T_1)\cap V(T_2)=S$. The basic problem of {\sc Steiner Tree
Packing} is defined as follows: the input consists of an undirected
graph $G$ and a subset of vertices $S\subseteq V(D)$, the goal is to
find a largest collection of edge-disjoint $S$-Steiner trees.  The Steiner tree packing problem has applications in VLSI circuit design \cite{Grotschel-Martin-Weismantel, Sherwani}. 
In this application, a Steiner tree is needed to share an electronic signal by a set of terminal nodes. Another application arises in the Internet Domain \cite{Li-Mao5}: Let a given graph $G$ represent a network, we choose arbitrary $k$ vertices as nodes such that one of them is a {\em broadcaster}, and all other nodes are either {\em users} or {\em routers} (also called {\em switches}). The broadcaster wants to broadcast as many streams of movies as possible, so that the users have
the maximum number of choices. Each stream of movie is broadcasted via a tree connecting all the users and the broadcaster. Hence we need to find the maximum number Steiner trees connecting all the users and the broadcaster, and it is a Steiner tree packing problem.

An {\em out-tree} (resp. {\em in-tree}) is an oriented tree in which every vertex
except one, called the {\em root}, has in-degree (resp.  out-degree) one.
An {\em out-branching} (resp. {\em in-branching}) of $D$ is a spanning out-tree
(resp.  in-tree) in $D$. For a digraph $D=(V(D), A(D))$, and a set $S\subseteq V(D)$ with $r\in S$ and $|S|\geq 2$, a {\em directed $(S, r)$-Steiner tree} or, simply, an {\em $(S, r)$-tree} is an out-tree $T$ rooted at $r$ with
$S\subseteq V(T)$ \cite{Cheriyan-Salavatipour}. Two $(S, r)$-trees are said to be {\em arc-disjoint} if they have no common arc. Two arc-disjoint $(S, r)$-trees are said to be {\em internally disjoint} if the set of common vertices of them is exactly $S$. Let $\kappa_{S,r}(D)$ (resp. $\lambda_{S,r}(D)$) be the maximum number of internally disjoint (resp. arc-disjoint) $(S, r)$-trees in $D$.

Cheriyan and Salavatipour \cite{Cheriyan-Salavatipour}, Sun and Yeo \cite{Sun-Yeo} studied the following two directed Steiner tree packing problems:

\noindent {\bf Arc-disjoint Directed Steiner Tree Packing (ADSTP):}
The input consists of a digraph $D$ and a subset of vertices
$S\subseteq V(D)$ with a root $r$, the goal is to find a largest
collection of arc-disjoint $(S, r)$-trees.

\noindent {\bf Internally-disjoint Directed Steiner Tree Packing (IDSTP):}
The input consists of a
digraph $D$ and a subset of vertices $S\subseteq V(D)$ with a root
$r$, the goal is to find a largest collection of internally disjoint
$(S, r)$-trees.

Now we introduce the concept of directed Steiner path packing which could be seen as a restriction of the directed Steiner tree packing problem. For a digraph $D=(V(D), A(D))$, and a set $S\subseteq V(D)$ with $r\in S$ and $|S|\geq 2$, a {\em directed $(S, r)$-Steiner path} or, simply, an {\em $(S, r)$-path} is a directed path $P$ started at $r$ with $S\subseteq V(P)$. Observe that the directed Steiner path is a generalization of the directed Hamiltonian path.

Two $(S, r)$-paths are said to be {\em arc-disjoint} if they have no common arc. Two arc-disjoint $(S, r)$-paths are said to be {\em internally disjoint} if the set of common vertices of them is exactly $S$. Let $\kappa^p_{S,r}(D)$ (resp. $\lambda^p_{S,r}(D)$) be the maximum number of internally disjoint (resp. arc-disjoint) $(S, r)$-paths in $D$. The directed Steiner path packing problems can be defined as follows:

\noindent {\bf Arc-disjoint Directed Steiner Path Packing (ADSPP):}
The input consists of a digraph $D$ and a subset of vertices
$S\subseteq V(D)$ with a root $r$, the goal is to find a largest
collection of arc-disjoint $(S, r)$-paths.

\noindent {\bf Internally-disjoint Directed Steiner Path Packing (IDSPP):}
The input consists of a
digraph $D$ and a subset of vertices $S\subseteq V(D)$ with a root
$r$, the goal is to find a largest collection of internally disjoint
$(S, r)$-paths.

The following concept of directed path connectivity is
related to directed Steiner path packing problem and is a natural
extension of path connectivity of undirected graphs (see \cite{Li-Mao5} for the introduction of path connectivity) to directed
graphs.
The {\em  directed path $k$-connectivity}
of $D$ is defined as
$$\kappa^p_k(D)= \min \{\kappa^p_{S,r}(D)\mid S\subseteq V(D), |S|=k, r\in
S\}.$$ Similarly, the {\em directed path $k$-arc-connectivity}
of $D$ is defined as
$$\lambda^p_k(D)= \min \{\lambda^p_{S,r}(D)\mid S\subseteq V(D), |S|=k, r\in
S\}.$$ By definition, when $k=2$, $\kappa^p_2(D)=\kappa(D)$ and
$\lambda^p_2(D)=\lambda(D)$. Hence, these two parameters could be seen
as generalizations of vertex-strong connectivity and arc-strong
connectivity of a digraph. The directed path $k$-connectivity and directed path $k$-arc-connectivity are also called {\em directed path connectivity}.

In Section~2, we study the complexity of the decision version of Steiner path packing problems in digraphs.
When both $k\geq 3, \ell\geq 2$ are fixed integers, we show that the problem of deciding whether there are at least $\ell$ internally disjoint directed $S$-Steiner paths in an Eulerian digraph $D$ is NP-complete, where $S\subseteq V(D)$ and $|S|=k$ (Theorem~\ref{thm1a}). However, when we consider the class of symmetric digraphs, the problem becomes polynomial-time solvable (Theorem~\ref{thm1c}). We also show that the problem of deciding whether there are at least $\ell$ arc-disjoint directed $S$-Steiner paths in a given digraph $D$ is NP-complete, where $S\subseteq V(D)$ and $|S|=k$ (Theorem~\ref{thm1b}).

In Section~3, we study the parameters $\kappa^p_k(D)$ and $\lambda^p_k(D)$. We show that the values $\lambda^p_k(D)$ is decreasing over $k$, but the values $\kappa^p_k(D)$ are neither increasing, nor decreasing over $k$. We give two upper bounds for $\kappa^p_k(D)$ in terms of classical connectivity $\kappa(D)$ and $\lambda(D)$  of digraphs (Theorem~\ref{thmb}). We also give sharp lower and upper bounds for the Nordhaus-Gaddum type relations of the parameter $\lambda^p_k(D)$ (Theorem~\ref{thmf}).

\section{Main results}

The problem of {\sc Directed $k$-Linkage} is defined as follows: Given a digraph $D$ and a (terminal) sequence $(s_1, t_1, \dots, s_k, t_k)$ of  distinct vertices of $D,$
decide whether $D$ has $k$ vertex-disjoint paths $P_1, \dots, P_k$, such that $P_i$ starts at $s_i$ and ends at $t_i$ for all $i\in [k].$

Sun and Yeo proved the NP-completeness of {\sc Directed 2-Linkage} for Eulerian digraphs.

\begin{thm}\label{thm101}\cite{Sun-Yeo}
The problem of {\sc Directed 2-Linkage} restricted to Eulerian digraphs is NP-complete.
\end{thm}



Now we can prove the NP-completeness of deciding whether $\kappa^p_{S,r}(D)\geq \ell$ for Eulerian digraphs (and therefore for general digraphs).

\begin{thm}\label{thm1a}
Let $k\geq 3, \ell \geq 2$ be fixed integers. For any Eulerian digraph $D$ and $S
\subseteq V(D)$ with $|S|=k$ and $r\in S$, the problem of deciding whether $\kappa^p_{S,r}(D) \geq \ell$ is NP-complete.
\end{thm}
\begin{pf}
It is not difficult to see that the problem belongs to NP. 
We will show that the NP-hardness of this problem  by reducing from the problem of {\sc Directed 2-Linkage} in Eulerian digraphs.
Let $[H; s_1,s_2,t_1,t_2]$ be an instance of {\sc Directed 2-Linkage} in Eulerian digraphs, that is, $H$ is an Eulerian digraph, and $(s_1, t_1, s_2, t_2)$ is a (terminal) sequence of  distinct vertices of $H$. 

We now produce a new Eulerian digraph $D$ as follows:

Let $V(H') = V(H) \cup S \cup \{r_1, r_2\},$
where $S=\{x_i \mid i\in [k]\}$, and let 
\[
\begin{array}{rcl}
 A(H') & = & A(H) \cup \{x_{1}s_1, t_1x_2, x_{k-1}s_2, t_2x_k, s_1r_1, r_1t_2, s_2r_2, r_2t_1\} \\
        &   & \cup \; \{ x_ix_{i+1} \mid i\in [k] \},\\
\end{array}
\]
where $x_{k+1}=x_1 $. Furthermore, we duplicate the arc $x_ix_{i+1}$ $\ell-1$ times for each $2\leq i\leq k-2$ or $i=k$, and duplicate the two arcs $x_1x_2$ and $x_{k-1}x_1$ $\ell-2$ times. Finally, to avoid parallel arcs, insert a new vertex $z_{i,i+1}^j$ to each arc of the form $x_ix_{i+1}$, where $j\in [\ell]$ if $2\leq i\leq k-2$ or $i=k$, and $j\in [\ell-1]$ otherwise. Let $D$ be the resulting digraph, and let $r=x_1$.

\begin{figure}[htb]
		\centering
		\begin{tikzpicture}
			\tikzset{arrow1/.style = {draw = black, thick, -{Latex[length = 3mm, width = 1.2mm]},}
			}
			
			\filldraw[black]    (0, 0)  circle (2pt)  node [anchor=north east] {$x_{1}$};
			\filldraw[black]    (1.5, 0)  circle (2pt)  node [anchor=south] {$z^{j}_{1,2}$};
			\filldraw[black]    (3, 0)  circle (2pt)  node [anchor=south] {$x_{2}$};
			\filldraw[black]    (5, 0)  circle (2pt)  node [anchor=south] {$x_{k-1}$};
			\filldraw[black]    (6.5, 0)  circle (2pt)  node [anchor=south] {$z^{j}_{k-1,k}$};
			\filldraw[black]    (8, 0)  circle (2pt)  node [anchor=north west] {$x_{k}$};
			\filldraw[black]    (4, 1.2)  circle (2pt)  node [anchor=south] {$z^{j}_{k,k+1}$};
			\filldraw[black]    (0.9, -2.5)  circle (2pt)  node [anchor=north east] {$s_{1}$};
			\filldraw[black]    (2.2, -2.5)  circle (2pt)  node [anchor=north east] {$t_{1}$};
			\filldraw[black]    (5.8, -2.5)  circle (2pt)  node [anchor=north east] {$s_{2}$};
			\filldraw[black]    (7.1, -2.5)  circle (2pt)  node [anchor=north west] {$t_{2}$};
			\filldraw[black]    (2, -5)  circle (2pt)  node [anchor=north east] {$r_{1}$};
			\filldraw[black]    (6.5, -5)  circle (3pt)  node [anchor=north west] {$r_{2}$};
			
			\draw[line width=1pt] (-0.4, -2) rectangle (8.4, -3.5);
			
			\draw[densely dotted] [line width=1.2pt]      (3.5, 0) -- (4.5, 0);
			
			\draw[arrow1] []      (0, 0) -- (1.5, 0);
			\draw[arrow1] []      (1.5, 0) -- (3, 0);
			\draw[arrow1] []      (5, 0) -- (6.5, 0);
			\draw[arrow1] []      (6.5, 0) -- (8, 0);
			\draw[arrow1] []      (0, 0) -- (0.9, -2.5);
			\draw[arrow1] []      (5, 0) -- (5.8, -2.5);
			\draw[arrow1] []      (2.2, -2.5) -- (3, 0);
			\draw[arrow1] []      (7.1, -2.5) -- (8, 0);
			\draw[arrow1] []      (0.9, -2.5) -- (2, -5);
			\draw[arrow1] []      (2, -5) -- (7.1, -2.5);
			\draw[arrow1] []      (6.5, -5) -- (2.2, -2.5);
			\draw[arrow1] []      (5.8, -2.5) -- (6.5, -5);
			\draw[arrow1] []      (8, 0) .. controls (6.5, 1) and (5, 1.2) .. (4, 1.2);
			\draw[arrow1] []      (4, 1.2) .. controls (3, 1.2) and (1.5, 1) .. (0, 0);
			
			\node at(4, -5.7) {$D$};
			\node at(8.7, -2.8) {$H$};
		\end{tikzpicture}
		\caption{The digraph $D$.}
		\label{figure1}
	\end{figure}
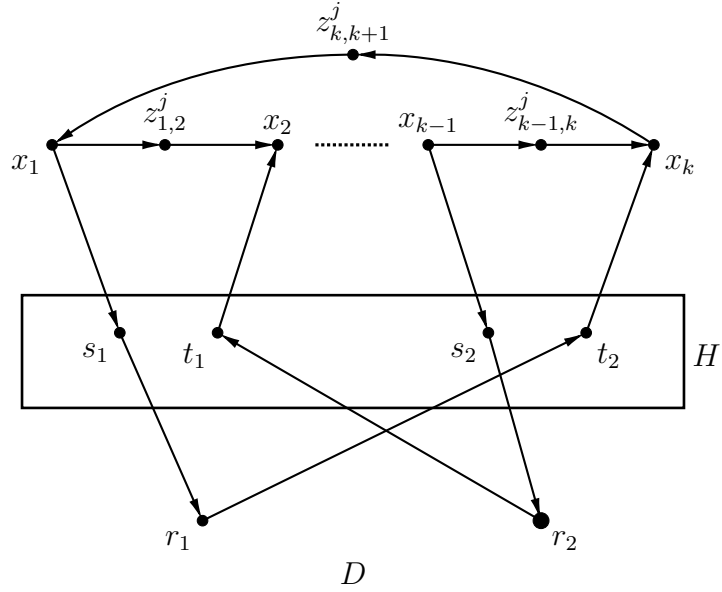

It can be checked that $D$ is Eulerian. We will show that $\kappa^p_{S,r}(D)\ge \ell$ if and only if there exist two vertex-disjoint paths, $Q_1$ and $Q_2$, in $H$ such that $Q_i$ is an $(s_i,t_i)$-path, for $i\in [2]$. This will complete the proof by Theorem~\ref{thm101}.

First assume that there exist two vertex-disjoint paths, $Q_1$ and $Q_2$ in $H$, such that $Q_i$ is an $(s_i,t_i)$-path, for $i\in [2]$.
Add the arcs $x_{1}s_1, t_1x_2, x_{k-1}s_2, t_2x_k$, and the path $x_2, z^{\ell}_{2,3}, x_3, \dots, z^{\ell}_{k-2, k-1}, x_{k-1}$ to $Q_1, Q_2$ and call the resulting path for $P_{\ell}$. Let $P_j$ be the path $x_1, z^{j}_{1,2}, x_2, \dots, z^{j}_{k-1, k}, x_{k}$, for each $j\in [\ell-1]$. It can be checked that the above $(S,r)$-paths $P_1, P_2,\ldots, P_{\ell}$ are internally disjoint, which implies that  $\kappa^p_{S,r}(D) \geq \ell$ as desired.

Conversely, assume that $\kappa^p_{S,r}(D) \geq \ell$, that is, there is a set of internally disjoint $(S, r)$-paths, say $\{P_j\mid j\in [\ell]\}$. Let $S'=S\cup \{z^j_{i,i+1}\mid j\in [\ell]~if~2\leq i\leq~k-2~or~i=k,~and~j\in [\ell-1]~otherwise\}$. By the construction of $D$, the subdigraph $D[S']$ is the union of the two paths $x_2, z^{\ell}_{2,3}, x_3, \dots, z^{\ell}_{k-2, k-1}, x_{k-1}$, and $x_k, z^{\ell}_{k,k+1}, x_1$, and $\ell-1$ arc-disjoint cycles: $x_1, z^{j}_{1,2}, x_2, \dots, z^{j}_{k-1, k}, x_{k}, z^{j}_{k, k+1}, x_1$ where $j\in [\ell-1]$. Since $deg_{D}^+(x_i)=\ell$ for $i\in [k-1]$ and $deg_{D}^{-}(x_i)=\ell$ for $2\leq i\leq k$, each $P_j$ contains precisely one out-neighbour of $x_i$ for $i\in [k-1]$ and one in-neighbour of $x_i$ for $2\leq i\leq k$. There must be two paths, say $Q_1, Q_2$, such that $Q_1$ contains $s_1, t_1$ (resp. $Q_2$ contains $s_2, t_2$), and starts from $x_1$ (resp. $x_{k-1}$) and ends at $x_2$ (resp. $x_{k}$) in $P_{j_1}$ (resp. $P_{j_2}$) for some $j_1\in [\ell]$ (resp. $j_2\in [\ell]$).

If $j_1=j_2$, then clearly there is a pair of vertex-disjoint $(s_1, t_1)$-path and $(s_2, t_2)$-path in $P_{j_1}$ (and therefore in $H$), as desired. It remains to consider the case that $j_1\neq j_2$. Without loss of generality, assume that $j_1=1, j_2=2$. Observe that $P_1=Q_1\cup Q'_1\cup Q''_1$ such that $Q'_1$ is an $(x_2, x_{k-1})$-subpath and $Q''_1$ is an $(x_{k-1}, x_{k})$-subpath in $P_1$, and $P_1=Q'_2\cup Q''_2\cup Q_2$ such that $Q'_2$ is an $(x_1, x_{2})$-subpath and $Q''_2$ is an $(x_{2},x_{k-1})$-subpath in $P_2$. Now let $P'_1=Q_1\cup Q'_1\cup Q_2$ and $P'_2=Q'_2\cup Q''_2\cup Q''_1$.
It can be checked that $\{P'_i, P_j\mid i\in [2], 3\leq j\leq \ell\}$ is still a set of internally disjoint $(S,r)$-paths. With a similar argument to that of the above case, we still have that
there is pair of vertex-disjoint $(s_1, t_1)$-path and $(s_2, t_2)$-path in $P'_1$ (and therefore in $H$).
\end{pf}



We can also prove the NP-completeness of deciding whether $\lambda^p_{S,r}(D)\geq \ell$ for general digraphs.

\begin{thm}\label{thm1b}
Let $k\geq 3, \ell \geq 2$ be fixed integers. For any digraph $D$ and $S
\subseteq V(D)$ with $|S|=k$ and $r\in S$, the problem of deciding whether $\lambda^p_{S,r}(D) \geq \ell$ is NP-complete.
\end{thm}
\begin{pf}
It is not difficult to see that the problem belongs to NP. 
We will show that the problem is NP-hard by reducing from {\sc Directed 2-Linkage} in Eulerian digraphs.
Let $[H; s_1,s_2,t_1,t_2]$ be an instance of {\sc Directed 2-Linkage} in Eulerian digraphs, that is, $H$ is an Eulerian digraph, and $(s_1, t_1, s_2, t_2)$ is a (terminal) sequence of  distinct vertices of $H$.

We first produce a digraph $D$ the same as that in Theorem~\ref{thm1a}. Secondly, we construct a new digraph $D'$ from $D$ as follows: replace every vertex $u$ of $H$ by two vertices $u^-$ and $u^+$ such that $u^-u^+$ is an arc in $D'$ and for every $uv\in A(H)$ add an arc $u^+v^-$ to $D'$. Furthermore, for $z\in S \cup \{r_1, r_2\}$, for every arc $zu$ in $D$ add an arc $zu^-$ to $D'$ and for every arc $uz$ add an arc $u^+z$ to $D'$.

It was proved in Theorem~\ref{thm1a} that $\kappa^p_{S,r}(D) \geq \ell$ if and only if $[H; s_1,s_2,t_1,t_2]$ is a positive instance of {\sc Directed 2-Linkage} in Eulerian digraphs. Therefore, to prove the theorem, it suffices to show that $\lambda^p_{S,r}(D') \geq \ell$ if and only if $\kappa^p_{S,r}(D) \geq \ell$.

We first assume that $\kappa^p_{S,r}(D) \geq \ell$, that is, there is a set of internally disjoint $(S,r)$-paths in $D$, say $\{P_i\mid i\in [\ell]\}$. By the construction of $D$, there must be two paths, say $Q_1, Q_2$, such that $Q_1$ contains $s_1, t_1$ (resp. $Q_2$ contains $s_2, t_2$), and starts from $x_{1}$ (resp. $x_{k-1}$) and ends at $x_2$ (resp. $x_{k}$) in $P_i$ (resp. $P_j$) for some $i\in [\ell]$ (resp. $j\in [\ell]$). Now we obtain a path $P'_i$ (resp. $P'_j$) in $D'$ from $P_i$ (resp. $P_j$) as follows: replace every vertex $u$ of $V(H)\cap V(P_i)$ (resp. $V(H)\cap V(P_j)$) by two vertices $u^-$ and $u^+$ such that $u^-u^+$ is an arc in $P'_i$ (resp. $P'_j$) and for every $uv\in A(H)\cap A(P_i)$ (resp. $uv\in A(H)\cap A(P_j)$) add an arc $u^+v^-$ to $A(P'_i)$ (resp. $A(P'_j)$). Furthermore, for each $z\in S \cup \{r_1, r_2\}$, for every arc $zu$ in $A(D)\cap A(P_i)$ (resp. $A(D)\cap A(P_j)$) add an arc $zu^-$ to $P'_i$ (resp. $P'_j$), and for every arc $uz$ add an arc $u^+z$ to $P'_i$ (resp. $P'_j$). Then combining with the remaining $(S,r)$-paths in $\{P_i\mid i\in [\ell]\}$, we obtain a set of arc-disjoint $(S,r)$-paths in $D'$, therefore $\lambda^p_{S,r}(D') \geq \ell$.

Conversely, we assume that $\lambda^p_{S,r}(D') \geq \ell$, that is, there is a set of arc-disjoint $(S,r)$-paths in $D'$, say $\{P_i\mid i\in [\ell]\}$. By the construction of $D'$, there must be two paths, say $Q'_1, Q'_2$, such that $Q'_1$ contains $s^-_1, s^+_1, t^-_1, t^+_1$ (resp. $Q'_2$ contains $s^-_2, s^+_2, t^-_2, t^+_2$), and starts from $x_{1}$ (resp. $x_{k-1}$) and ends at $x_2$ (resp. $x_{k}$) in $P_i$ (resp. $P_j$) for some $i\in [\ell]$ (resp. $j\in [\ell]$). Observe that the remaining paths in $\{P_i\mid i\in [\ell]\}$ only use arcs/vertices from $D[S']$, where $S'$ is defined in the proof of Theorem~\ref{thm1a}.
Also, observe that $Q'_1$ must be of the form $x_{1}, s^-_1, s^+_1, a^-_1, a^+_1, \dots, a^-_p, a^+_p, t^-_1, t^+_1, x_2$ and $Q'_2$ must be of the form $x_{k-1}, s^-_2, s^+_2, b^-_1, b^+_1, \dots, b^-_q, b^+_q, t^-_2, t^+_2, x_{k}$ and furthermore, we have $\{a_i\mid i\in [p]\}\cap \{b_j\mid j\in [q]\}=\emptyset$. Now we obtain $P'_i$ (resp. $P'_j$) from $P_i$ (resp. $P_j$) by replacing $Q'_1$ (resp. $Q'_2$) with the path $x_{1}, s_1, a_1, \dots, a_p, t_1, x_2$ (resp. $x_{k-1}, s_2, b_1, \dots, b_p, t_2, x_{k}$). It can be checked that, combining with the remaining $(S,r)$-paths in $\{P_i\mid i\in [\ell]\}$, we get a set of $\ell$ internally disjoint $(S,r)$-paths in $D$, therefore $\kappa^p_{S,r}(D) \geq \ell$.
\end{pf}

Recall that in Theorem~\ref{thm1a}, we showed that when $D$ is an Eulerian digraph, the problem of deciding whether $\kappa^p_{S,r}(D)\geq \ell$ with $|S|=k$ is NP-complete, where both $k\geq 3, \ell\geq 2$ are fixed integers. However, when we consider the class of symmetric digraphs, the problem becomes polynomial-time solvable. We need the following result by Sun and Yeo.

\begin{lem}\label{thmsym}\cite{Sun-Yeo}
Let $D$ be a symmetric digraph and let $s_1,s_2,\ldots,s_r, t_1,t_2,\ldots,t_r$ be vertices in $D$ (not necessarily disjoint) and
let $S \subseteq V(D)$. We can in $O(|V(G)|^3)$ time decide if there for all $i=1,2,\ldots,r$ exists an $(s_i,t_i)$-path, $P_i$,
such that no internal vertex of any $P_i$ belongs to $S$ or to any path $P_j$ with $j \not=i$ (the end-points of $P_j$ can also not
be internal vertices of $P_i$).
\end{lem}

By Lemma~\ref{thmsym}, we obtain the following lemma which will be used in the argument of Theorem~\ref{thm1c}.

\begin{lem}\label{thmc1}
Let $k\geq 3$ and $\ell \geq 2$ be fixed integers and let $D$ be a symmetric digraph. Let $S \subseteq V(D)$ with $|S|=k$, and let $r$ be an arbitrary vertex in $S$. Let $\{A_i\mid 0\leq i\leq \ell\}$ be a partition of $A[S]$, where $A[S]$ denotes the arc set in $D[S]$.

We can in time $O(n^3k^{k\ell})$ decide if there exists a set of $\ell$ internally disjoint $(S,r)$-paths, say $\{P_i \mid i\in [{\ell}]\}$, such that $A(P_i) \cap A[S] = A_i$ for each $i\in [\ell]$ (note that $A_0$ are the arcs in $A[S]$ not used in any of the paths).
\end{lem}
\begin{pf}
Let $P$ be any $(S,r)$-path in $D$, we define the {\em skeleton} of $P$ as the path we obtain from $P$ by contracting all vertices in $V(P)\setminus S$.
Let $P^s$ be a skeleton of an $(S,r)$-path in $D$. Note that $V(P^s)=S$ and therefore there are at most $k^k$ different skeletons of $(S,r)$-paths in $D$.

Our algorithm will try all possible $\ell$-tuples, ${\cal P}^s = (P_1^s, P_2^s, \ldots, P_{\ell}^s)$, of skeletons of $(S,r)$-paths and determine if there is a set of $\ell$ internally disjoint $(S,r)$-paths, say $\{P_i \mid i\in [{\ell}]\}$, such that $A(P_i) \cap A[S] = A_i$ and  $P_i^s$ is the skeleton of $P_i$ for each $i\in [\ell]$. If such a set of paths exists for any ${\cal P}^s$, then we return this solution, and if no such set of paths exist for any ${\cal P}^s$, then we return that no solution exists. We will prove that this algorithm gives the correct answer and compute its time complexity.

If our algorithm returns a solution, then clearly a solution exists. So now assume that a solution exists and let $\{P_i \mid i\in [{\ell}]\}$ be the desired set of internally disjoint $(S,r)$-paths. When we consider ${\cal P}^s = (P_1^s, P_2^s, \ldots, P_{\ell}^s)$, where $P_i^s$ is the skeleton of $P_i$, our algorithm will find a
solution, so the algorithm always returns a solution if one exists.


Given such an $\ell$-tuples, ${\cal P}^s$, we need to determine if there is a set of $\ell$ internally disjoint $(S,r)$-paths, say $\{P_i \mid i\in [{\ell}]\}$, such that $A(P_i) \cap A[S] = A_i$ and  $P_i^s$ is the skeleton of $P_i$ for each $i\in [\ell]$.
We first check that the arcs in $A_i$ belong to the skeleton $P_i^s$ and that no vertex in $V(D)\setminus S$ belongs to more than one skeleton.
If the above does not hold, then the desired paths do not exist, so we assume that the above holds in the following argument. For every arc $uv \not\in A[S]$ that belongs to some skeleton $P_i^s$, we want to find a $(u,v)$-path in $D - A[S]$,
such that no internal vertex on any path belongs to $S$ or to a different path.
This can be done in $O(n^3)$ time by Lemma~\ref{thmsym}.  If such paths exist, then we obtain the desired $(S,r)$-paths by substituting each $uv$ by the corresponding $(u,v)$-path. Otherwise, the desired set of $(S,r)$-paths does not exist. 

Therefore, the algorithm works correctly and has complexity $O(n^3k^{k\ell})$ since the number of different $\ell$-tuples, ${\cal P}^s$, that we need to consider is bounded by the function $(k^k)^{\ell}=k^{k\ell}$.
\end{pf}

We will now prove the polynomiality for $\kappa^p_{S,r}(D)$ on symmetric digraphs.

\begin{thm}\label{thm1c}
Let $k\geq 3$ and $\ell \geq 2$ be fixed integers. We can in polynomial time decide if $\kappa^p_{S,r}(D) \geq \ell$ for any symmetric digraph $D$ with $S \subseteq V(D)$, with $|S|=k$ and $r\in S$.
\end{thm}
\begin{pf}
Let $D$ be a symmetric digraph with $S \subseteq V(D)$, where $|S|=k$. Let  ${\cal A} = \{A_i\mid 0\leq i\leq \ell\}$ be a partition of $A[S]$. By Lemma~\ref{thmc1}, we can in time $O(n^3k^{k\ell})$ decide if there exists a set of $\ell$ internally disjoint $(S,r)$-paths, say $\{P_i \mid i\in [{\ell}]\}$, such that $A(P_i) \cap A[S] = A_i$ for each $i\in [\ell]$.

We will now use the algorithm of Lemma~\ref{thmc1} for all possible partitions ${\cal A} = \{A_i\mid 0\leq i\leq \ell\}$. If we find the desired set of $\ell$ internally disjoint $(S,r)$-paths for any such a partition, then we return ``$\kappa^p_{S,r}(D) \geq \ell$''. Otherwise, we return ``$\kappa^p_{S,r}(D) < \ell$''.
Note that if  $\kappa^p_{S,r}(D) \geq \ell$, then we will correctly determine that $\kappa^p_{S,r}(D) \geq \ell$, when we consider the correct partition ${\cal A}$, which proves that the above algorithms will always return the correct answer.

Hence, we deduce that the algorithm works correctly and has complexity $O(n^3k^{k\ell}(\ell+1)^{k^2/2})$ which is a polynomial in $n$, since the number of partitions ${\cal A}$ of $A(D[S])$ is bounded by $(\ell+1)^{|A(D[S])|} \leq (\ell+1)^{k^2/2}$, and the two parameters $k$ and $\ell$ are fixed.
\end{pf}

\section{The parameters $\kappa^p_k(D)$ and $\lambda^p_k(D)$}

The following proposition can be verified using definitions of
$\kappa^p_k(D)$ and $\lambda^p_{k}(D)$.

\begin{pro}\label{proa}
Let $D$ be a digraph of order $n$, and let $k$ be an integer such that $2 \leq k \leq n$.
Then the following assertions hold:
\begin{description}
\item[(1):] $\lambda^p_{k+1}(D)\leq \lambda^p_{k}(D)$ when $k \leq n-1$.

\item[(2):] $\kappa^p_k(D')\leq \kappa^p_k(D)$ and $\lambda_k(D')\leq \lambda_k(D)$
where $D'$ is a spanning subdigraph of $D$.

\item[(3):] $\kappa^p_k(D) \leq \lambda^p_k(D) \leq \min\{\delta^+(D), \delta^-(D)\}$.
\end{description}
\end{pro}

The following example means that Proposition~\ref{proa}(1) (the monotonicity of $\lambda^p_k(D)$ over $k$) doesnot hold for the parameter $\kappa^p_k(D)$.

\vspace{2mm}

\noindent{\bf Example~1:} We now define a digraph $D$ on eight vertices, such that $\kappa^p_8(D)=\kappa^p_2(D)=2$, but $\kappa^p_4(D)\leq 1$, which shows that
the values $\kappa^p_k(D)$ are neither increasing, nor decreasing over $k$. 

Let $V(D)=X \cup Y \cup Z$, where $X=\{x_1,x_2\}$, $Y=\{y_1,y_2,y_3\}$
and $Z=\{z_1,z_2,z_3\}$. Furthermore let $D$ contain all arcs from $Z$ to $X$, all arcs from $Z$ to $Y$ and all arcs from $X$ to $Y$. Finally add all arcs $Y$ to $Z$, except $y_iz_i$ for $i=1,2,3$. 

It is not difficult to show that $\kappa^p_2(D)=\kappa(D)=2$.
We now prove that $\kappa^p_{V(D),r}\geq 2$ for any $r\in V(D)$. When $r\in \{x_1, x_2\}$, say $r=x_1$, let $P_1: x_1y_1z_2y_2z_3y_3z_1x_2$ and $P_2: x_1y_2z_1y_3z_2y_1z_3x_2$. When $r\in \{y_1, y_2, y_3\}$, say $r=y_1$, let $P_1: y_1z_2x_1y_2z_3x_2y_3z_1$ and $P_2: y_1z_3x_1y_3z_2x_2y_2z_1$. When $r\in \{z_1, z_2, z_3\}$, say $r=z_1$, let $P_1: z_1x_1y_1z_2x_2y_2z_3y_3$ and $P_2: z_1x_2y_1z_3x_1y_3z_2y_2$. Observe that in each case, $P_1$ and $P_2$ are arc-disjoint Hamiltonian paths started at $r$. Hence, $\kappa^p_{V(D),r}\geq 2$ for any $r\in V(D)$, which means that $\kappa^p_8(D)\geq 2$. Furthermore, by Proposition~\ref{proa}(3), $\kappa^p_8(D)\leq \delta^-(D)=2$. Hence, $\kappa^p_8(D)=2$.

Let $S=Z \cup \{x_1\}$ and let $r=x_1$. For the sake of contradiction assume that there exist two internally disjoint $(S,r)$-paths, $P_1$ and $P_2$. Either $|V(P_1) \cap Y| \leq 1$ or $|V(P_2) \cap Y| \leq 1$, as $|Y|=3$. Without loss of generality assume that $|V(P_1) \cap Y| \leq 1$, which implies that there is a vertex in $Z = S \setminus \{r\}$ with no arc into it in $P_1$, a contradiction. Therefore $\kappa^p_4(D) \leq 1$, as desired.\qed

\vspace{2mm}

We still need the following famous result on the Hamiltonian decomposition of complete digraphs.

\begin{thm}(Tillson's decomposition theorem)\cite{Tillson}\label{Tillson}
The arcs of $\overleftrightarrow{K}_n$ can be decomposed into
Hamiltonian cycles if and only if $n\neq 4,6$.
\end{thm}

By Theorem~\ref{Tillson}, $\overleftrightarrow{K}_n$ is the union of $n-1$ arc-disjoint Hamiltonian cycles if $n\neq 4,6$, and so by Proposition~\ref{proa}(3), we directly have that $$\lambda^p_k(\overleftrightarrow{K}_n)=\kappa^p_n(\overleftrightarrow{K}_n)=n-1$$ for any $2\leq k\leq n$, when $n\geq 7$. The means that the bound in Theorem~\ref{thmb}(ii) is sharp when $k=n$.

\begin{thm}\label{thmb}
Let $2\leq k\leq n$ be an integer. The following assertions hold:
\begin{description}
\item[(i)] $\kappa^p_k(D)\leq \kappa(D)$ when $n\ge \kappa(D)+k$.
\item[(ii)] $\kappa^p_k(D) \leq \lambda(D)$.
\end{description}
\end{thm}
\begin{pf}
\noindent{\bf Part (i)}.
For $k=2$, we have $\kappa^p_2(D)= \kappa(D)$ by definition. In the following argument we therefore consider the case of $k\ge 3$.
If $\kappa(D)=0$, then $D$ is not strong and $\kappa^p_k(D)=0$, as can be seen by letting $r,x \in S$ be chosen such that there is
no $(r,x)$-path in $D$.
If $\kappa(D)=n-1$, then we have $\kappa^p_k(D)\leq n-1=\min\{\delta^+(D), \delta^-(D)\}$ by
Proposition~\ref{proa}(3), so we may assume that $1 \leq \kappa(D) \leq n-2$. There now exists a $\kappa(D)$-vertex
cut, say $Q$, for two vertices $u,v$ in $D$ such that there is no
$(u,v)$-path in $D-Q$. Let $S=\{u,v\}\cup S'$ where $S'\subseteq
V(D)\setminus (Q\cup \{u,v\})$ and $|S'|=k-2$. Observe that in each
$(S, u)$-path, the $u-v$ path must contain a vertex in $Q$. By the
definition of $\kappa^p_{S,r}(D)$ and $\kappa^p_k(D)$, we have
$\kappa^p_k(D)\leq \kappa^p_{S,r}(D)\leq |Q|=\kappa(D)$. 

\vspace{2mm}

\noindent{\bf Part (ii)}. Let $D$ be any digraph with $\ell =\lambda(D)$.
Let $A' \subseteq A(D)$ be defined such that $|A'|=\ell$ and $D-A'$ is not strong.
Let  $r$ and $x$ be vertices in $D$ such that there is no $(r,x)$-path in $D-A'$.
Let  $S \subseteq V(D)$ be chosen such that  $\{r,x\} \subseteq S$ and $|S|=k$.
Observe that in each $(S,r)$-path in $D$ we must use an arc from $A'$, which implies that
$\kappa^p_k(D) \leq \kappa^p_{S,r}(D) \leq |A'| = \lambda(D)$.  
\end{pf}

Given a graph parameter $f(G)$, the Nordhaus-Gaddum Problem is to
determine sharp bounds for (1) $f(G) + f(G^c)$ and (2) $f(G)f(G^c)$.
The Nordhaus-Gaddum type
relations have received wide attention; see a survey paper \cite{Aouchiche-Hansen} by Aouchiche and Hansen. The following Theorem~\ref{thmf}
concerns such type of a problem for the parameter $\lambda^p_k(D)$.

\begin{thm}\label{thmf}
For a digraph $D$ with order $n\geq 7$, the following assertions hold:
\begin{description}
\item[(i)]~$0\leq \lambda^p_k(D)+\lambda^p_k(D^c)\leq n-1$. Moreover, both bounds are sharp. 
\item[(ii)]~$0\leq \lambda^p_k(D){\lambda^p_k(D^c)}\leq
\lfloor(\frac{n-1}{2})^2\rfloor$. Moreover, both bounds are sharp.
\end{description}
\end{thm}
\begin{pf}
\noindent{\bf Part (i)}. Since $D\cup D^c=\overleftrightarrow{K}_n$,
Proposition~\ref{proa}(3) implies the following, 
\[
\lambda^p_k(D)+\lambda^p_k(D^c)\leq   \delta^+(D) + \delta^+(D^c) \leq d_D^+(x) + d_{D^c}^+(x) = n-1,
\]
where $x \in V(D)$ is arbitrary.

If $H\cong \overleftrightarrow{K}_n~(n\geq 7)$, then we have $\lambda^p_k(H)=n-1$ and
$\lambda^p_k(H^c)=0$, so the upper bound is sharp. The lower bound is
clear. Furthermore, the lower bound holds, if and only if $\lambda^p_k(D)=\lambda^p_k(D^c)=0$. As an example, consider a non-strong tournament, $T$, in which case it is not difficult to see that $T^c$ is also non-strong, so now $\lambda^p_k(T)=\lambda^p_k(T^c)=0$ by Proposition~\ref{proa}(1).

\vspace{2mm}

\noindent{\bf Part (ii)}. The lower bound is clear. Furthermore, the lower bound holds, if and only if $\lambda^p_k(D)=0$ or
$\lambda^p_k(D^c)=0$. For the sharpness, we can just use the example in the argument of part~(i).

For the upper bound, we have

\[
\lambda^p_k(D){\lambda^p_k(D^c)}\leq
\left(\frac{\lambda^p_k(D)+\lambda^p_k(D^c)}{2}\right)^2\leq
\left(\frac{n-1}{2}\right)^2.
\]

Since both $\lambda^p_k(D)$ and $\lambda^p_k(D^c)$ are integers, the upper bound holds.
We now prove the sharpness of the upper bound. By Theorem~\ref{Tillson}, the complete digraph $\overleftrightarrow{K}_n$ can be decomposed into $n-1$ arc-disjoint Hamilton cycles, when $n$ is odd. Let $D$ consist of $(n-1)/2$ of these arc-disjoint Hamilton cycles, which implies that $\lambda^p_k(D)=\lambda^p_k(D^c)=(n-1)/2$, which shows the sharpness of the upper bound. This completes the proof.
\end{pf}

\vskip 1cm

\noindent{\bf Acknowledgement.} This work was supported by Zhejiang Provincial Natural Science Foundation of China under Grant No. LY23A010011 and Yongjiang Talent Introduction Programme of Ningbo under Grant No. 2021B-011-G.

\end{document}